\documentclass[a4paper,11pt]{amsart}
\usepackage{amsfonts}
\usepackage{amsmath}
\usepackage{amssymb}
\usepackage{amsthm}
\usepackage{url}
\usepackage[english]{babel}
\usepackage{verbatim}
\usepackage{graphicx}
\usepackage{array}
\usepackage{blkarray}
\usepackage[caption=false]{subfig}
\usepackage{multirow}
\usepackage{enumitem, hyperref}

\makeatletter
\def\namedlabel#1#2{\begingroup
    #2%
    \def\@currentlabel{#2}%
    \phantomsection\label{#1}\endgroup
}
\makeatother

\numberwithin{equation}{section}

\usepackage{graphicx}
\DeclareGraphicsExtensions{.jpg,.gif,.eps}
\def\qed{\hfill $\diamondsuit$}

\theoremstyle{plain}
\newtheorem{theorem}{Theorem}[section]

\theoremstyle{remark}
\newtheorem{remark}{Remark}
\theoremstyle{definition}

\newtheorem{assumption}{Assumption}

\DeclareMathOperator*{\argmaxC}{arg\,max}

\begin{document}

   \author{Jun Maeda}
   \address{University of Warwick, Coventry, United Kingdom}
   \email{J.Maeda@warwick.ac.uk}
   \author{Saul D. Jacka}
   \thanks{Saul Jacka gratefully acknowledges funding received from the EPSRC grant EP/P00377X/1 and is also grateful to the Alan Turing Institute for their financial support under the EPSRC 
grant EP/N510129/1.}
   \address{University of Warwick, Coventry, United Kingdom\\
   and
The Alan Turing Institute, London,  United Kingdom }
   \email{S.D.Jacka@warwick.ac.uk}

  \title[Evaluation of the rate of convergence in the PIA]{\textit{Evaluation of the rate of convergence in the PIA}}

   \begin{abstract}
 Folklore says that Howard's Policy Improvement Algorithm converges extraordinarily fast, even for controlled diffusion settings.
 
In a previous paper, we proved that approximations of the solution of a particular parabolic partial differential equation obtained via the policy improvement algorithm  show a quadratic local convergence. 

In this paper, we show that we obtain the same rate of convergence of the algorithm in a more general setup. This provides some explanation as to why the  algorithm converges fast. 

We provide an example by solving a semilinear elliptic partial differential equation numerically by applying the algorithm and check how the approximations  converge to the analytic solution. 
\end{abstract}



   \date{\today}


\maketitle

\section{Introduction}
\label{section: Introduction}
In \cite{1.}, we introduced a new model for pricing derivatives products when we have a position concentration in the over-the-counter market. The model requires us to solve a nonlinear partial differential equation (PDE) and this may prevent traders from using it in practice due to possible difficulties in implementing a solution in their pricing models. To overcome this difficulty, we used the policy improvement algorithm (PIA) to enable us to approximate the nonlinear PDE by a series of linear ones parameterized by a control. The solutions of the linear PDEs converge to that of the original semilinear PDE as we iteratively solve the linear PDEs under the algorithm. Since their stochastic volatility pricing models can solve linear PDEs (as in Heston's model), the traders can now implement the new model. We further showed that the PIA approximated solutions show \textit{quadratic local convergence} (QLC) to the analytic solution. This provides an explanation of why the convergence happens so fast.

The natural question to ask is how general this QLC is in the PIA framework. In this paper, we consider a general infinite time horizon problem and calculate the rate of convergence of the PIA-derived approximations to that of the corresponding semilinear elliptic PDE. We give three conditions which enable us to show the QLC. These assumptions are indeed satisfied by the problem considered in \cite{1.}.  We describe in Remark~\ref{remark: 1} how some of these assumptions can be relaxed.
 
The rest of the paper is organized as follows:  Section~\ref{section: Setup} briefly explains the setup. In Section~\ref{section: Test_Case}, we state the main theorem about  Quadratic  Local Convergence of the  approximated solutions to the semilinear PDE. We give a numerical example in Section~\ref{section: Example} and give some concluding questions in Section~\ref{section: Conclusion}.

\section{Setup}
\label{section: Setup}
We briefly explain our setup. 

Let $(\Omega, \mathcal{F}, (\mathcal{F}_t)_{t\ge0}, P)$ be a filtered probability space. We assume that $\mathcal{E}$ is a simply connected, convex, and bounded subset of ${\mathbb{R}}^n$ that has $C^{2, \beta}$ boundary. 
We define

\begin{equation}
\tau_\mathcal{E}(\mathcal{Y}) :=\inf\{t\geq 0; \mathcal{Y}_t \notin \mathcal{E}\}
\end{equation}

for any continuous process $\mathcal{Y} = (\mathcal{Y}_t)_{t\geq0}$.

For a control $\Pi$ and starting point $z$, we wish to define the controlled process $Z^{z, \Pi}$ by
\begin{equation}\label{eq: SDE_0}
Z^{z, \Pi}_t = z + \int^t_0 \sigma(Z^{z, \Pi}_s, \Pi_s)dB_s + \int^t_0 \mu(Z^{z, \Pi}_s, \Pi_s)ds \quad 0\leq t \leq \tau_\mathcal{E}(Z^{z, \Pi})\text{,}
\end{equation}
where  $\sigma: {\mathbb{R}}^n\times {\mathbb{R}}^d \rightarrow \mathbb{R}^{n\times n}$ and $\mu: {\mathbb{R}}^n \times {\mathbb{R}}^d \rightarrow {\mathbb{R}}^n$ are measurable mappings, $B$ is an $n$-dimensional Wiener process and $\Pi$ takes values in $A=\mathbb{R}^d$.

For any $z \in {\mathbb{R}}^n$ define $\mathcal{A}(z)$, the set of admissible control at $z$, as
\begin{align}
\begin{split}
\mathcal{A}(z) &:= \{\Pi = (\Pi_t)_{t\ge0}; \Pi \text{ is adapted to } (\mathcal{F}_t)_{t\ge0}\text{,} \Pi_t(\omega) \in {\mathbb{R}}^d\text{ for every } t\geq 0 \text{ and } \omega \in \Omega\text{,}\\
&\text{and there exists a process } Z^{z, \Pi}=(Z^{z, \Pi}_t)_{t\ge0}\text{ that satisfies } \eqref{eq: SDE_0} \text{ and is unique in law}\}\text{.}
\end{split}
\end{align}

A measurable function $\pi: \Omega \times (0, \infty] \rightarrow {\mathbb{R}}^d$ is a \textit{Markov policy} if for every $z \in \Omega$ and $\forall T>0$ there exists a process $Z^{z, \pi}_t$ that is unique in law and satisfies the following:

\begin{align}\label{eq: SDE}
\begin{split}
Z^{z, \pi}_t &= z + \int^t_0 \sigma(Z^{z, \pi}_s, \pi(Z^{z, \pi}_s, s))dB_s + \int^t_0 \mu(Z^{z, \pi}_s, \pi(Z^{z, \pi}_s, s))ds\\
&= z + \int^t_0 \sigma_\pi(Z^{z, \pi}_s)dB_s + \int^t_0 \mu_\pi(Z^{z, \pi}_s)ds \quad 0 \leq t \leq T\wedge \tau_\mathcal{E}\text{.}
\end{split}
\end{align}

We define the payoff function $V^\Pi$  for any admissible $\Pi$ as

\begin{align}\label{eq: v_pi_definit}
\begin{split}
V^\Pi(z) :&= E\bigg(\int^{\tau_\mathcal{E}}_0 e^{-\alpha{t}}f(Z^{z, \Pi}_t, \Pi_t)dt + e^{-\alpha(\tau_\mathcal{E})}g(Z^{z, \Pi}_{\tau_\mathcal{E}})\bigg)\\
&= E\bigg(\int^{\tau_\mathcal{E}}_0 e^{-\alpha{t}}f^{\Pi_t}(Z^{z, \Pi}_t)dt+ e^{-\alpha(\tau_\mathcal{E})}g(Z^{z, \Pi}_{\tau_\mathcal{E}})\bigg)\text{,}
\end{split}
\end{align}

where $\alpha$ is some positive constant and $f: {\mathbb{R}}^n \times {\mathbb{R}}^d \rightarrow {\mathbb{R}}$ and $g: {\mathbb{R}}^n \rightarrow \mathbb{R}$. We assume that $f^\pi$ is $C^2$ with respect to $\pi$ and $g$ is continuous. The problem is to find the value function $V$ defined as

\begin{equation}
V := \sup_{\Pi \in \mathcal{A}}V^\Pi \text{.}
\end{equation}

For any Markov policy $\pi$ that is Lipschitz continuous on $\bar{\mathcal{E}}$, define $L^\pi:C^2\rightarrow C$ by 
\begin{equation}\label{eq: ell_def}
L^\pi \phi := \frac{1}{2}Tr\{\sigma_\pi^T (H \phi)\sigma_\pi\} + \mu_\pi^T\nabla\phi = \sum_{i,j}a^\pi_{ij}\frac{\partial^2 \phi}{\partial x_i \partial x_j} + \sum_i b_i \frac{\partial \phi}{\partial x_i},
\end{equation}
where $H\phi$ is the Hessian of $\phi$.

From \cite{2.}, $V^\pi$ satisfies the PDE

\begin{equation}
L^\pi V^\pi - \alpha V^\pi + f^\pi = 0\text{.}
\end{equation}

Starting from a Markov policy $\pi_0$, the PIA defines successive controls by the recursion 
\begin{equation}
\pi_{i+1} = \displaystyle\argmaxC_{a \in {A}} \big(L^{a}V^{\pi_i} - \alpha V^{\pi_i} + f^a\big) \text{.}
\end{equation}

Note finally that, we assume that  $\exists \nu>0$ such that the differential operator $L^\pi$ is uniformly elliptic, i.e.,

\begin{equation}\label{eq: ellipticity_condition}
\frac{1}{\nu} |\xi|^2 \le a^\pi_{ij}\xi_i\xi_j \le \nu |\xi|^2 \qquad \forall \xi_i, \xi_j \in \mathcal{E}, \forall \pi\in \mathbb{R}^d \text{.}
\end{equation}

\section{Main Results}
\label{section: Test_Case}
We make the following assumptions in this section.

\begin{assumption}\label{assumption: 1}
$\mu_\pi$ is in the form of ${\mathbf{M}}\pi + {\mathbf{b}}$ for some constant $n\times d$ matrix $\mathbf{M}$ and $n$ dimensional vector $\mathbf{b}$.
\end{assumption}

\begin{assumption}\label{assumption: 2}
$\sigma_\pi$ is independent of $\pi$.
\end{assumption}

\begin{assumption}\label{assumption: 3}
$f^\pi$ is strictly and uniformly concave in $\pi$, i.e. $\exists \lambda>0$ such that ${\bf{x}}^T(H_\pi f^\pi) {\bf{x}} \le -\lambda ||{\bf{x}}||^2<0$ for all ${\bf{x}}\in  \bar{\mathcal{E}}$, where $H_\pi f^\pi$ represents the Hessian of $f^\pi$ with respect to $\pi$.
\end{assumption}

With these assumptions, we show the following:

\begin{theorem}\label{prop: one}
Under Assumptions~\ref{assumption: 1}, \ref{assumption: 2}, and \ref{assumption: 3}, there exists a $C>0$ such that 

\begin{equation}\label{eq: prop_formula}
\lVert V^{\pi_{i+1}} - V^{\pi_i} \rVert_{2, \beta} \le C \lVert V^{\pi_i} - V^{\pi_{i-1}}\rVert_{2, \beta}^2\text{,}
\end{equation}

where $C$ only depends on the domain $\mathcal{E}$, the ellipticity constant $\nu$ from \eqref{eq: ellipticity_condition}, and the bounds on the coefficients of the differential operator $L^\pi$ (see \cite{3.} for the definition of the Sobolev norm $||\cdot||_{2,\beta}$).
\end{theorem}
\begin{remark}\label{q con}
Applying \eqref{eq: prop_formula} iteratively, we obtain
\begin{equation}\label{eq: iteratition_inequality}
\lVert V^{\pi_{i+1}} - V^{\pi_i} \rVert_{2, \beta} \le \frac{\{{C}\lVert V^{\pi_{1}} - V^{\pi_0} \rVert_{2, \beta}\}^{2^i}}{C}\text{.}
\end{equation}
Therefore, once $\lVert V^{\pi_{i+1}} - V^{\pi_i} \rVert_{2, \beta}<1/{C}$, convergence is {\em extremely} fast.
\end{remark}
\begin{remark}\label{remark: 1}
Suppose that ${A}$,  the action space (the value space for $\pi$) is not $\mathbb{R}^d$. We may replace it by its image under $\mu(x,\cdot)$, provided we simultaneously replace $f$ by $\tilde f$ given by
$$\tilde f(x, m)=\sup_{\pi \in (x,\cdot)^{-1}(m)}f(x, \pi),$$
since we wish to maximise $V^\pi$.
Suppose that this image is $\mathcal{M}$. By allowing relaxed controls (see for example \cite{dje}) we can replace this by $\mathcal{N}:=\bar{co}(\mathcal{M})$, the closure of the convex hull of $\mathcal{M}$.
This will simultaneously replace  $\tilde f(x,\cdot)$ by $\bar{f}$, the smallest concave majorant of $\tilde f$. 
If $\tilde f$ is strictly uniformly concave and $\mathcal{N}$ is an affine set in $\mathbb{R}^n$ then we recover Assumptions 1 and 3.

As we shall see, the proof of Theorem~\ref{prop: one} relies heavily on Taylor's theorem and the disappearance of $\nabla^a(L^a V^{\pi_i}-\alpha V^{\pi_i}+f^a)$ at its maximum. 
So, if $\mathcal{N}$ is a compact subset of ${\mathbb{R}}^d$  then we hit a problem when the maximizer $\mu$ lies on the boundary of $\mathcal{N}$.

 We should still be able to obtain good approximations to $V$ with QLC  by extending the action space to $\tilde{\mathcal{N}}$ and extending $\bar f^\mu$ to ${f}^{*,\mu}$ in such a way that $\bar f^\mu=f^{*,\mu}$ in ${\mathcal{N}}$, ${f}^{*,\mu}$ always takes its maximum, $\hat f$ in the interior of $\tilde{\mathcal{N}}$ and $\hat f-\sup_{\mu\in \mathcal{N}}\bar f^\mu\leq \epsilon$.
\end{remark}

\begin{remark}
We considered the elliptic case, but the parabolic case follows exactly in the same fashion. We have the following theorem:

\begin{theorem}\label{thm: 2}
Under the same assumptions as in Theorem~\ref{prop: one} with a given initial condition, \eqref{eq: prop_formula} holds in the parabolic case.
\end{theorem}

This is a generalization of Proposition 5.3 in \cite{1.}.
\end{remark}

The proof of Theorem~\ref{prop: one} is deferred to the Appendix.

\section{Numerical Example}
\label{section: Example}

We apply the PIA in solving numerically a semilinear elliptic PDE.

We take $\mathcal{E} \subset \mathbb{R}^2$ to be $[0.5,2.0]\times[0.5,2.0]$ with its corners smoothed in a $C^2$ fashion (this is needed to apply the boundary estimate in Theorem \ref{prop: one}). The SDEs we consider are

\begin{equation}\label{eq: SDE_example}
	\left\{
	\begin{array}{c}
	dx =  \pi{x}\,dt + \sigma{x}\,dW^1\text{,}\\
	dy = \pi{y} \,dt +  \eta{y}\,dW^2\text{,}\\
	<dW^1,dW^2> = 0 \text{,}
\end{array}
\right.
\end{equation}

where $W^1$ and $W^2$ are 1-dimensional Wiener processes and $\pi \in \mathbb{R}$.

Thus

\begin{equation}
\mu_\pi = 
\begin{pmatrix}
\pi{x}\\
\pi{y}
\end{pmatrix}
\quad \text{ and } \quad
\sigma_\pi = 
\begin{pmatrix}
\sigma x && 0\\
0 && \eta{y}\\
\end{pmatrix}
\text{.}
\end{equation}

We take $f^\pi$ to be

\begin{equation}
f^\pi = 1 - \frac{1}{2}\pi^2\text{.}
\end{equation}

We define $V^{\pi}$ as in \eqref{eq: v_pi_definit} with $g \equiv 0$ on $\partial \mathcal{E}$.

Then, $V^{\pi_{i}}$ satisfies the elliptic PDE:

\begin{align}\label{eq: linearPDE}
\begin{split}
\frac{1}{2}\sigma^2 x^2\frac{\partial^2 V^{\pi_{i}}}{\partial x^2} + \frac{1}{2}\eta^2 y^2\frac{\partial^2 V^{\pi_{i}}}{\partial y^2} +  \pi_i {x}\frac{\partial V^{\pi_{i}}}{\partial x} + \pi_i {y} \frac{\partial V^{\pi_{i}}}{\partial y} - \alpha{V^{\pi_{i}}} +1 - \frac{1}{2}\pi_i^2 = 0\text{,}
\end{split}
\end{align}

where $\pi_i$ is determined by

\begin{equation}\label{eq: numerical_example_Pidef}
\pi_i = {x}\frac{\partial V^{\pi_{i-1}}}{\partial x} +  {y} \frac{\partial V^{\pi_{i-1}}}{\partial y}\text{.}
\end{equation}

Note that if $V^{\pi_i}$ converges, the limit function $V$ satisfies a semilinear elliptic PDE

\begin{align}\label{eq: semilinearPDE}
\begin{split}
\frac{1}{2}\sigma^2 x^2\frac{\partial^2 V}{\partial x^2} + \frac{1}{2}\eta^2 y^2\frac{\partial^2 V}{\partial y^2} - \alpha{V}+1 - \frac{1}{2}\bigg(x\frac{\partial V}{\partial x} + y\frac{\partial V}{\partial y}\bigg)^2= 0\text{.}
\end{split}
\end{align}

The variables we use are in Table~\ref{table: parameters}. 

\begin{table}[ht]
\caption{Parameters we use for the numerical calculation.}
\label{table: parameters}
\begin{center}
    \begin{tabular}{ |c| c|}
    \hline
    parameter & value\\ 
    \hline
    $\alpha$ & 0.03 \\ 
    $\sigma$ & 2.0\\
    $\eta$ & 0.2 \\ 
    $x_{max}$ & 2.0 \\ 
    $x_{min}$ & 0.50 \\ 
    $y_{max}$ & 2.0 \\  
    $y_{min}$ & 0.50 \\ 
    ToleranceLevel1& 0.00001 \\ 
    ToleranceLevel2& 0.001 \\ 
    discretization nodes& 100 \\
    \hline
    \end{tabular}
\end{center}
\end{table}

We use the explicit finite difference method (FDM) to see the convergence starting at $\pi_0 \equiv 0$ with the boundary condition $V^{\pi_i} |_{\partial \mathcal{E}}\equiv 0$. We discretize \eqref{eq: linearPDE} and obtain

\begin{align}\label{eq: PIA_difference}
\begin{split}
&\frac{1}{2}\bigg(\frac{\sigma^2 x_{j}^2}{\Delta{x}^2}+\frac{\pi_i(j, k) x_j}{\Delta{x}}\bigg)V(j+1, k) + \frac{1}{2}\bigg(\frac{\sigma^2 x_{j}^2}{\Delta{x}^2}-\frac{\pi_i(j, k) x_j}{\Delta{x}}\bigg)V(j-1, k)\\&+\frac{1}{2}\bigg(\frac{\eta^2 y_{k}^2}{\Delta{y}^2} + \frac{\pi_i(j, k) y_k}{\Delta{y}}\bigg)V(j, k+1)+\frac{1}{2}\bigg(\frac{\eta^2 y_{k}^2}{\Delta{y}^2} - \frac{\pi_i(j, k) y_k}{\Delta{y}}\bigg)V(j, k-1)\\
&-\bigg(\frac{\sigma^2 x_j^2}{\Delta{x}^2} + \frac{\eta^2 y_k^2}{\Delta{y}^2} + \alpha\bigg)V(j, k)+1 -\frac{1}{2}\pi_i^2(j, k)\\
&=p_{j+1, k}V(j+1, k) + p_{j-1, k}V(j-1, k) + p_{j, k+1}V(j, k+1) + p_{j, k-1}V(j, k-1)\\
&+p_{j,k} V(j, k) + q_i(j,k)= 0\text{,}
\end{split}
\end{align}

where $x_j$ and $y_j$ represent coordinates of the mesh points, $V(j, k)$ and $\pi_i(j, k)$ are corresponding values at the mesh points, and $\Delta{x}$ and $\Delta{y}$ are corresponding mesh size. We therefore can write \eqref{eq: PIA_difference} in the form

\begin{align}\label{eq: PIA_eq_system}
\begin{split}
V(j, k) = -\frac{1}{p_{j, k}}\big\{&p_{j+1, k}V(j+1, k) + p_{j-1, k}V(j-1, k) \\
&+ p_{j, k+1}V(j, k+1) + p_{j, k-1}V(j, k-1) + q_i(j, k)\big\}\text{.}
\end{split}
\end{align}

We use the Gauss-Seidel method \cite{4.} together with the PIA to solve \eqref{eq: semilinearPDE}. 
The procedure is as follows:

\begin{enumerate}
\item
Set $V^{0}(j, k)=0$ and $\pi_0(j, k)=0 \quad \forall (j, k)$.
\item
Assume that we have $\pi_i$ and $V^{\ell}$ for all the mesh points. Use \eqref{eq: PIA_eq_system} to calculate the values $V^{\ell+1}(j, k)$. That is, use
\begin{align}\label{eq: PIA_eq_system_elltoell1}
\begin{split}
V^{\ell + 1}(j, k) = -\frac{1}{p_{j, k}}\big\{&p_{j+1, k}V^{\ell}(j+1, k) + p_{j-1, k}V^{\ell}(j-1, k) \\
&+ p_{j, k+1}V^{\ell}(j, k+1) + p_{j, k-1}V^{\ell}(j, k-1) + q_i(j, k)\big\}
\end{split}
\end{align}
to calculate $V^{\ell+1}(j, k)$.
\item
Iteratively solve for $V^{\ell + 1}$ from $V^{\ell}$ and stop when $\displaystyle\max_{j,k}|V^{\ell+1}(j,k) - V^{\ell}(j,k)| < \text{ToleranceLevel1}$. Calculate $\pi_{i+1}(j,k)$ by
\begin{equation}\label{eq: numerical_example_Pidef_itoi+1}
\pi_{i+1}(j,k) = \frac{V^{\ell +1}(j+1, k) - V^{\ell+1}(j-1, k)}{2\Delta x}x_j + \frac{V^{\ell +1}(j, k+1) - V^{\ell+1}(j, k-1)}{2\Delta y}y_k\text{.}
\end{equation}
\item
Repeat Procedure (3) and end the program when $\displaystyle\max_{j,k}|\pi_i(j,k) - \pi_i(j,k)| < \text{ToleranceLevel2}$. The numerical solution to  \eqref{eq: semilinearPDE} is $V^{\ell+1}(j, k)$.
\end{enumerate}

The method converges  if  the diagonal terms of the matrix are greater than the sum of the absolute values of the off-diagonal terms (Theorem 4.4.5, \cite{5.}). That is, on \eqref{eq: PIA_eq_system}, the method converges if

\begin{equation}\label{eq: GS_sufficient_cond}
|p_{j, k}| > |p_{j+1, k}| + |p_{j-1, k}|+ |p_{j, k+1}| + |p_{j, k-1}|\text{.}
\end{equation}

With $\pi_i$ small enough, the condition of the cited theorem is satisfied with the parameters we have chosen.

To compare the calculation load, we also numerically solved the corresponding linear PDE

\begin{equation}\label{eq: linearPDE_example}
\frac{1}{2}yx^2\frac{\partial^2 V}{\partial x^2} + \frac{1}{2}y\eta^2\frac{\partial^2 V}{\partial y^2} - \alpha{V} - 1= 0\text{.}
\end{equation}

The only difference between \eqref{eq: semilinearPDE} and \eqref{eq: linearPDE_example} is the existence of the term $-(1/2)\{x({\partial V}/{\partial x}) + y({\partial V}/{\partial y})\}^2$.

Table~\ref{table: comparison1} shows the numerical results in both linear and semilinear cases. For the linear case \eqref{eq: linearPDE_example}, we used Gauss-Seidel method with the tolerance level equal to ToleranceLevel1 in Table~\ref{table: parameters}. We see that the linear and semilinear cases have similar order in terms of the number of calculations to approximate to the specified tolerance level.
\begin{table}[ht]
\caption{Calculation load comparison for successful convergence. 1 calculation here means solving the difference equation \eqref{eq: PIA_eq_system_elltoell1} once at one point.}
\label{table: comparison1}
\vspace{.1 in}
\begin{center}
    \begin{tabular}{|c|c|c|}
    \hline
    Problem Type & Method & \# of calculations \\ 
    \hline
    Linear & FDM (Gauss-Seidel) & 24,541,704 \\ 
    Semilinear & PIA \& Gauss-Seidel & 34,372,107 \\ 
    \hline
    \end{tabular}
\end{center}
\end{table}

Table~\ref{table: pia_detail} shows the result in more detail.

\begin{table}[ht]
\caption{Detail of the calculations in the PIA.}
\label{table: pia_detail}
\vspace{.1 in}
\begin{center}
    \begin{tabular}{|c|c|c|c|c|}
    \hline
    PIA  & Max Difference in& Max Difference in& \# of  & Calculation\\ 
    steps & $|\pi_{i} - \pi_{i-1}|$ & $|V^{\pi_{i}}- V^{\pi_{i-1}}|$&calculations & time\\ 
    \hline
    0 & 2.15455038 & & 24,541,704 & 0:16\\ 
    1 & 1.55932909 &0.02563695 & 8,017,218 & 0:05 \\ 
    2 & 0.16986263 &0.00372773 & 1,744,578 & 0:01 \\
    3 & 0.00400477 &0.00006031 & 58,806 & 0:00 \\ 
    4 & 0.00066038 & 0.00000995& 9,801 & 0:00 \\ 
    \hline
    \end{tabular}
\end{center}
\end{table}

\begin{figure}[htp]
	\centering
	\subfloat[]{%
		\includegraphics[width=0.33\linewidth, height =0.3\textheight]{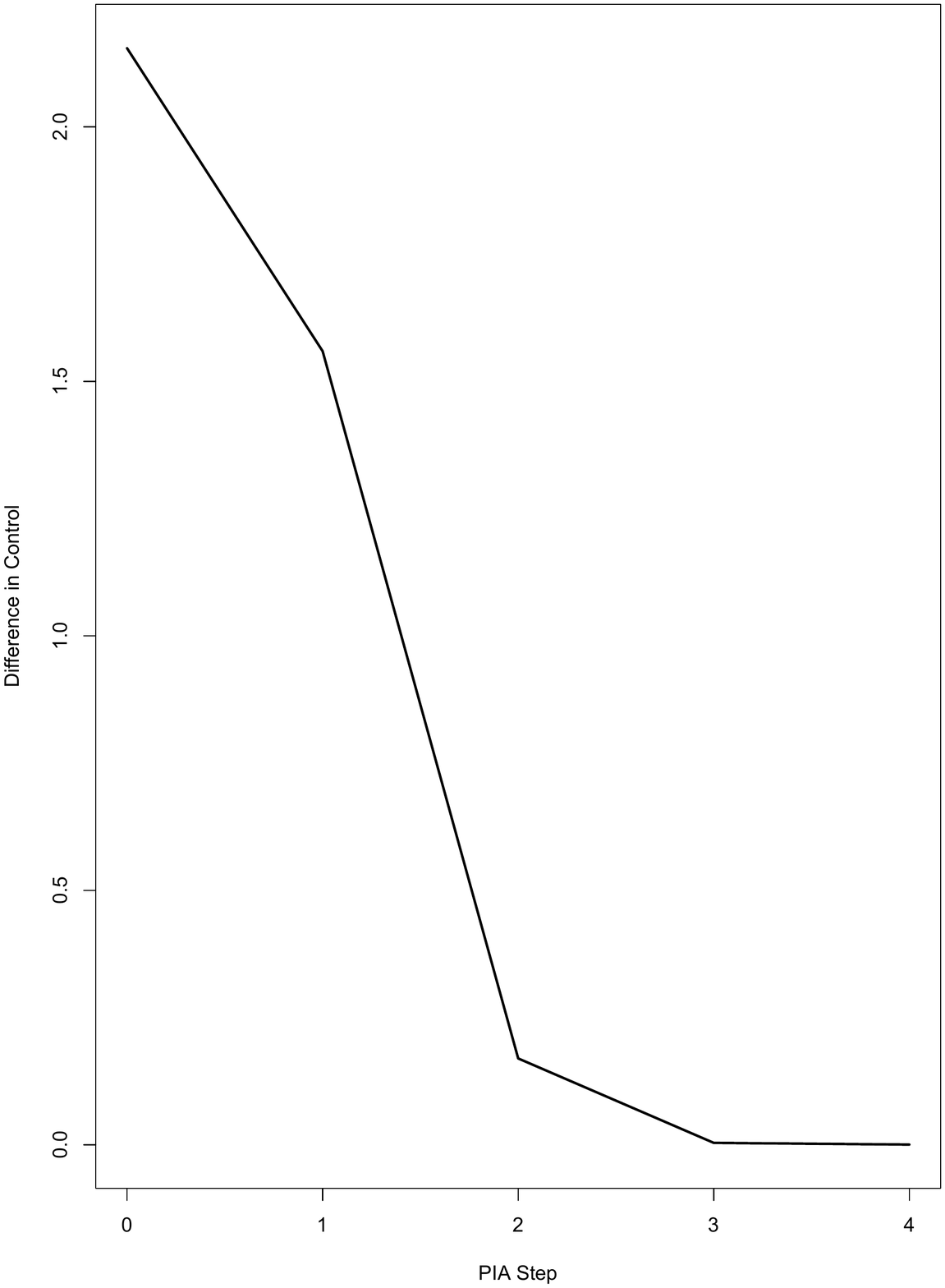}%
		\label{fig: Fig1.eps}%
		}%
	\hfill
	\subfloat[]{%
		\includegraphics[width=0.33\linewidth, height =0.3\textheight]{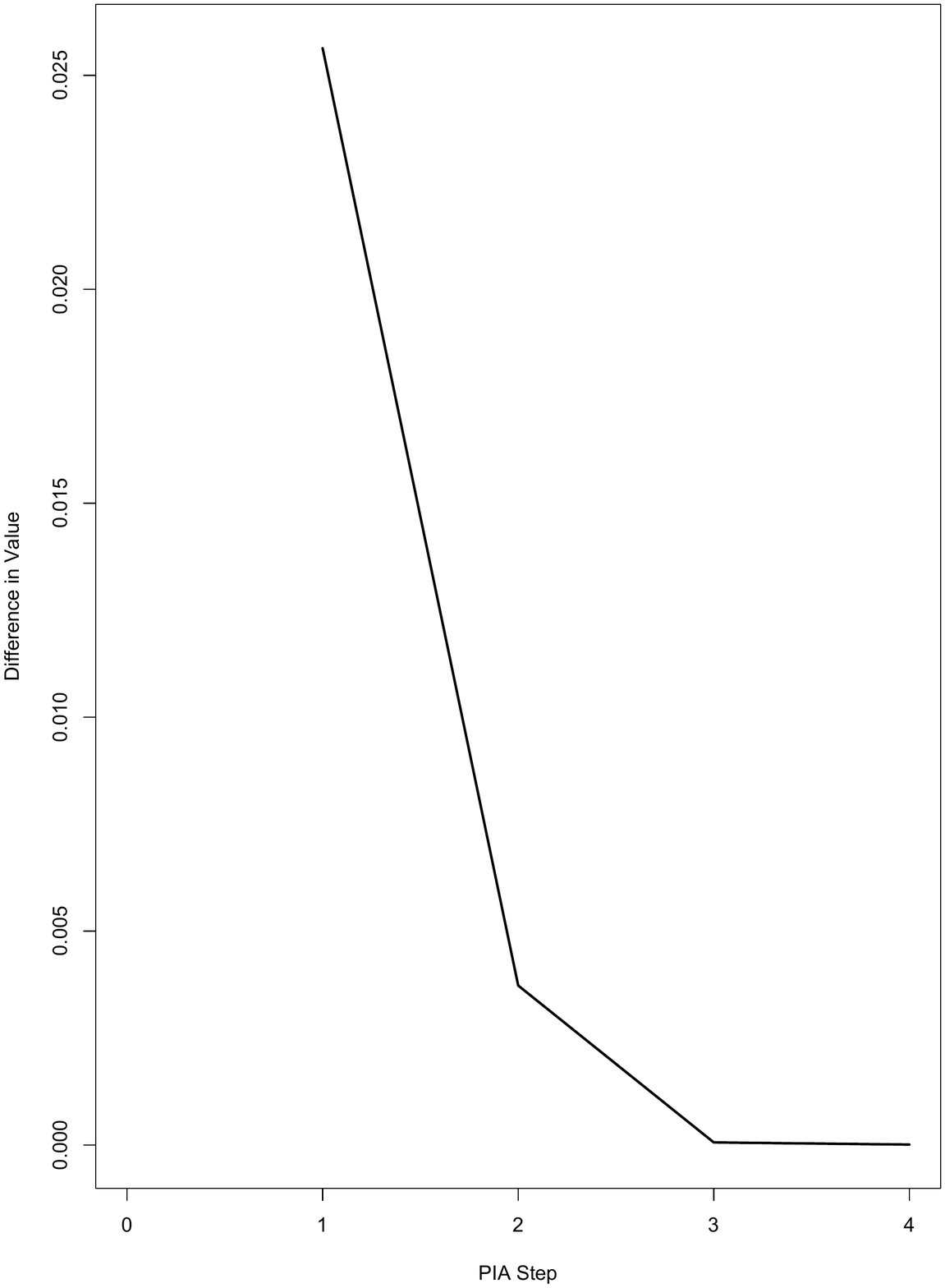}%
		\label{fig: Fig2.eps}%
		}%
	\hfill
	\subfloat[]{%
		\includegraphics[width=0.33\linewidth, height =0.3\textheight]{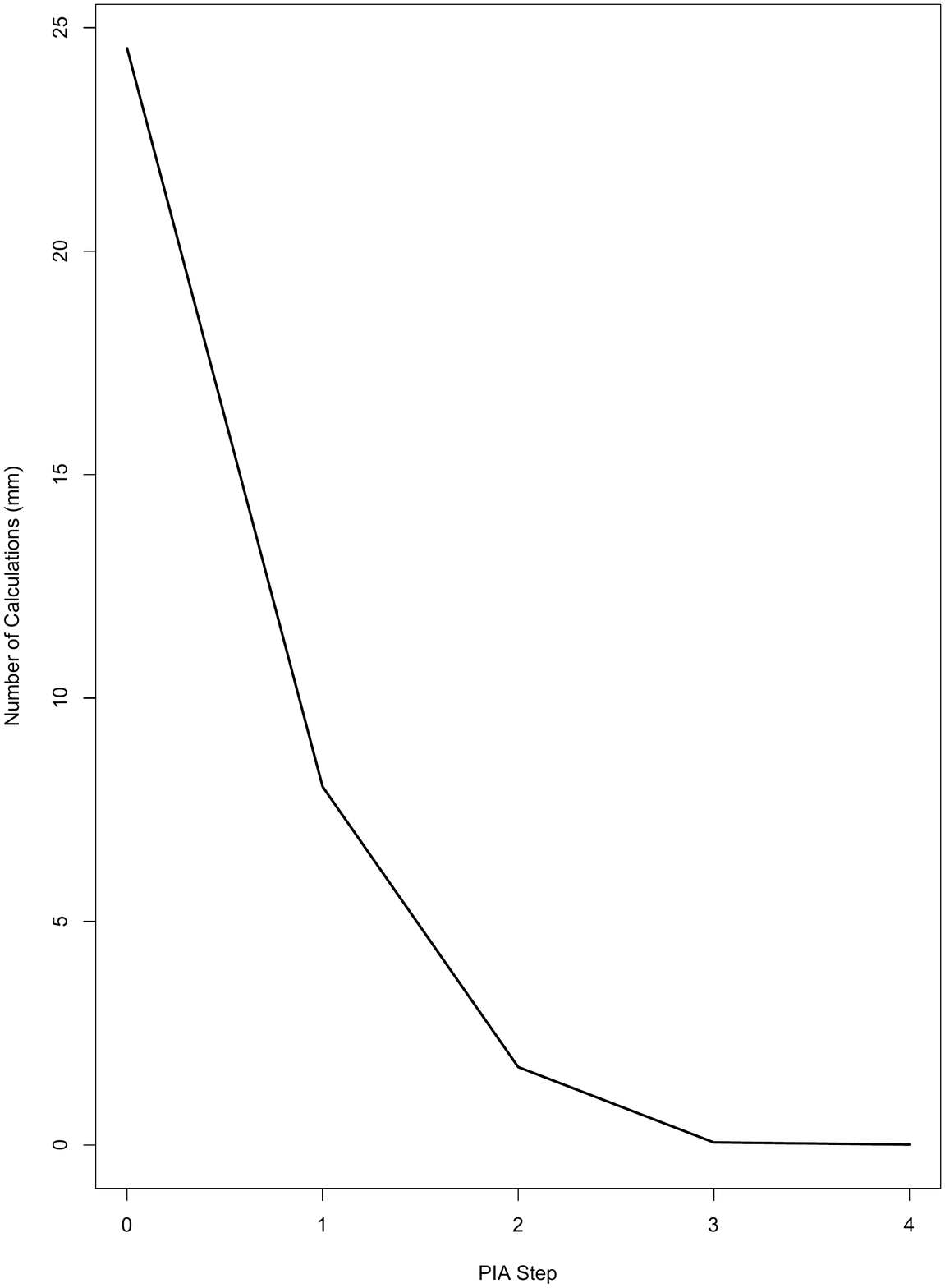}%
		\label{fig: Fig3.eps}%
		}%
	\caption{Graphs of the data in Table~\ref{table: pia_detail}. (A) the maximum of  $|\pi_{i} - \pi_{i-1}|$ in each step, (B) the maximum of $|V^{\pi_{i}} - V^{\pi_{i-1}}|$ in each step, and (C) the number of calculations in each step.}
	\label{fig: semilinear.png}
\end{figure}

Table~\ref{table: pia_detail} shows that the first step in the PIA already decreases the number of calculation to get the convergence in the Gauss-Seidel method. The data is plotted in Figure~\ref{fig: semilinear.png}.

\begin{remark}
We did try applying the FDM directly to the differential equation \eqref{eq: semilinearPDE}, but could not get the convergence.
\end{remark}

\section{Conclusions}
\label{section: Conclusion}

We have shown that the PIA has the QLC property in a fairly general framework. The natural questions to ask are
\begin{itemize}
\item [1.] Can we show  QLC under weaker conditions?

and 
\item[2.] Can we show some convergence rate outside the ``local quadratic region'' (see Remark \ref{q con})? 
\end{itemize}




\renewcommand{\thesubsection}{\Alph{subsection}}
\subsection*{Appendix}
\numberwithin{figure}{subsection}
\numberwithin{equation}{subsection}
\numberwithin{theorem}{subsection}

\subsection{Proof of Theorem~\ref{prop: one}}\label{App: Appendix}
$V^{\pi_i}$ satisfies 

\begin{equation}\label{eq: PDE_prop}
\frac{1}{2}Tr\{\sigma^T (H V^{\pi_i})\sigma\} + \mu_{\pi_i}^T\cdot\nabla V^{\pi_i} -\alpha V^{\pi_i} + f^{\pi_i}=0,
\end{equation}

and since Assumption~\ref{assumption: 2} is that $\sigma$ does not depend on $\pi$, $\pi_i$ is determined by the iteration:

\begin{align}\label{eq: defPi}
\begin{split}
\pi_{i+1} &= \argmaxC_{\pi \in \mathcal{A}} \bigg(\frac{1}{2}Tr\{\sigma^T(H {V^{\pi_{i}}}) \sigma\}+ \mu_{\pi}^T \cdot \nabla{V^{\pi_{i}}} + f^{\pi}\bigg)\\
& = \argmaxC_{\pi \in \mathcal{A}} \bigg(\mu_{\pi}^T\cdot\nabla{V^{\pi_{i}}} + f^{\pi}\bigg)\text{.}
\end{split}
\end{align}

From Assumption~\ref{assumption: 1}, we can write

\begin{equation}\label{eq: explicit_mu}
\mu_\pi = \mathbf{M}\pi + \mathbf{b} \text{.}
\end{equation} 

It then follows from \eqref{eq: defPi} that

\begin{equation}\label{eq: argmax_derivative}
\mathbf{M}^T\nabla V^{\pi_n} + \nabla_\pi f^\pi|_{\pi = \pi_{n+1}} = {\bf{0}} \text{.}
\end{equation}

Subtracting \eqref{eq: argmax_derivative} with $n=i-1$ from the same equation with $n=i$, and setting $W_i := V^{\pi_{i+1}}-V^{\pi_{i}}$, we obtain

\begin{equation}\label{eq: argmax_derivative2}
\mathbf{M}^T \nabla W_{i-1} + \nabla_\pi f^\pi|_{\pi = \pi_{i+1}} - \nabla_\pi f^\pi|_{\pi = \pi_{i}}= {\bf{0}}\text{.}
\end{equation}

Using the Mean Value Theorem, we can then write \eqref{eq: argmax_derivative2} as

\begin{equation}\label{eq: argmax_derivative3}
\mathbf{M}^T\nabla W_{i-1} + (H_\pi f^\pi)|_{\pi'}^T\cdot (\pi_{i+1} - \pi_i)= \bf{0}
\end{equation}

for some $\pi'\in {\mathbb{R}}^d$.

It follows from Assumption~\ref{assumption: 3} that  $H_\pi f^\pi$ is negative definite, hence invertible, so we can rewrite \eqref{eq: argmax_derivative3} as

\begin{equation}\label{eq: pi_increment}
\pi_{i+1} - \pi_i = - \{(H_\pi f^\pi)|_{\pi'}^T\}^{-1}\mathbf{M}^T\nabla W_{i-1}\text{.}
\end{equation}

Comparing \eqref{eq: PDE_prop} for $i$ and $i+1$,

\begin{equation}\label{eq: iterative_pde}
	\left\{
	\begin{array}{c}
	\frac{1}{2}Tr\{\sigma^T (H V^{\pi_{i+1}})\sigma\} + (\mathbf{M}\pi_{i+1}+ \mathbf{b})^T\cdot\nabla V^{\pi_{i+1}}-\alpha V^{\pi_{i+1}} + f^{\pi_{i+1}}=0 \text{,}\\
	\frac{1}{2}Tr\{\sigma^T (H V^{\pi_{i}})\sigma\} + (\mathbf{M}\pi_{i} + \mathbf{b})^T\cdot \nabla V^{\pi_{i}} -\alpha V^{\pi_{i}} + f^{\pi_{i}}=0\text{,}
\end{array}
\right.
\end{equation}

and subtracting, we get

\begin{align}\label{eq: PDE_W}
\begin{split}
\frac{1}{2}Tr\{\sigma^T (H W_i)\sigma\}  & + (\mathbf{M}\pi_{i+1}  + \mathbf{b})^T\cdot \nabla W_i -\alpha W_i \\
&+ \{\mathbf{M}(\pi_{i+1} - \pi_i)\}^T\cdot \nabla V^{\pi_{i}}+ (f^{\pi_{i+1}} - f^{\pi_{i}})=0\text{.}
\end{split}
\end{align}

We define $\mathcal{R}_i$ as

\begin{equation}
\mathcal{R}_i = (\pi_{i+1} - \pi_i)^T\mathbf{M}\nabla V^{\pi_{i}}+ (f^{\pi_{i+1}} - f^{\pi_{i}})\text{,}
\end{equation}
then we obtain, from Taylor's theorem,

\begin{align}\label{eq: residual_R}
\begin{split}
\mathcal{R}_i &=  (\pi_{i+1} - \pi_i)^T\cdot\bigg\{ \mathbf{M}\nabla V^{\pi_{i}}+ \nabla_\pi f^\pi|_{\pi = \pi_{i+1}} -\frac{1}{2}H_\pi f^\pi|_{\pi'}\cdot(\pi_{i+1} - \pi_i)\bigg\}\\
& =  -\frac{1}{2}(\pi_{i+1} - \pi_i)^T(H_\pi f^\pi)|_{\pi'}^T\cdot(\pi_{i+1} - \pi_i)\text{ from \eqref{eq: argmax_derivative}.}
\end{split}
\end{align}

Using \eqref{eq: pi_increment}, we can write \eqref{eq: residual_R} as

\begin{equation}\label{eq: residual_R_2}
\mathcal{R}_i =  -\frac{1}{2}(\mathbf{M}^T \nabla W_{i-1})^T(H_\pi f^\pi|_{\pi'})^{-1}(\mathbf{M}^T\nabla W_{i-1})\text{,}
\end{equation}
and then we can rewrite \eqref{eq: PDE_W} as

\begin{equation}\label{eq: PDE_W2}
\frac{1}{2}Tr\{\sigma^T (H W_i)\sigma\}  + (\mathbf{M}\pi_{i+1}+ \mathbf{b})^T\nabla W_i -\alpha W_i + \mathcal{R}_i=0\text{.}
\end{equation}

We have the same Dirichlet condition on the boundary of the domain for each $V^\pi$, therefore $W_i \equiv 0$ on $\partial \mathcal{E}$. From Schauder's estimate on second order linear elliptic partial differential equations \cite[pg. 108]{3.}, we conclude that
\begin{equation}\label{eq: Schauder_est}
\lVert W_i \rVert_{2, \beta}  \le C \lVert \mathcal{R}_i \rVert_{0, \beta}= C\lVert \nabla W_{i-1}\rVert_{0, \beta}^2 \le C \lVert W_{i-1}\rVert_{2, \beta}^2\text{,}
\end{equation}

where the constant $C$ depends only on the domain $\mathcal{E}$, the ellipticity constant $\nu$, and the bounds on the coefficients of the elliptic differential operator.
\qed

\end{document}